\documentclass{svmult}
\usepackage{mathptmx}
\usepackage{helvet}
\usepackage{courier}
\usepackage{type1cm}

\usepackage{latexsym}
\usepackage{amsmath}
\usepackage{amsfonts}
\usepackage{mathrsfs}

\newcommand{\bz}{ {\bf 0} }

\newcommand{\diag}{\mbox{diag} }
\newcommand{\rank}{\mbox{rank} }

\newcommand{\Rs}{\mathbb{R}}

\newcommand{\Sn}{{\cal S}_n }
\newcommand{\snn}{{\cal S}_{n-1} }

\newcommand{\Sh}{{\cal S}_H }

\newcommand{\mm}{\bar{m} }

\newcommand{\rr}{\bar{r} }

\newcommand{\yy}{\hat{y} }

\newcommand{\KK}{{\cal K} }

\newcommand{\TT}{{\cal T} }
\newcommand{\LL}{{\cal L} }

\newcommand{\EE}{{\cal E} }

\newcommand{\M}{{\cal M} }
\newcommand{\trace}{{\rm trace\,}}

\newcommand{\bpr}{{\bf Proof.} \hspace{1 em}}
\newcommand{\epr}{ \\ \hspace*{4.5in} $\Box$ }
\newcommand{\beq}{ \begin{equation} }
\newcommand{\eeq}{ \end{equation} }
\newcommand{\bt}{ \begin{tabular} }
\newcommand{\et}{ \end{tabular} }

\begin{document}

\title*{Universal Rigidity of Bar Frameworks in General Position:}
\subtitle{A Euclidean Distance Matrix Approach}

\author{ A. Y. Alfakih}

\institute{Department of Mathematics and Statistics,
          University of Windsor,
          Windsor, Ontario N9B 3P4,
          Canada. \email{alfakih@uwindsor.ca }}

\maketitle



\abstract{
A configuration $p$ in $r$-dimensional Euclidean space is a finite
collection of labeled points $p^1,\ldots, p^n$ in $\Rs^r$ that affinely
span $\Rs^r$. Each configuration $p$ defines a Euclidean distance matrix
$D_p =(d_{ij})$ = $(||p^i-p^j||^2)$, where
$||\; ||$ denotes the Euclidean norm.
A fundamental problem in distance geometry is to find out
whether or not, a given proper subset of the entries of $D_p$
suffices to uniquely determine the entire matrix $D_p$.
This problem is known as the universal rigidity problem of bar frameworks.
In this chapter, we present a unified approach for the universal rigidity
of bar frameworks, based on Euclidean distance matrices (EDMs), or equivalently,
on projected Gram matrices. This approach makes the universal rigidity problem
amenable to semidefinite programming methodology. Using this approach, we
survey some recently obtained results and their proofs, emphasizing
the case where the points $p^1,\ldots, p^n$ are in general position.}

\section{Introduction}

A {\em configuration} $p$ in $r$-dimensional Euclidean space is a
finite collection of labeled points $p^1,\ldots,p^n$ in $\Rs^r$
that affinely span $\Rs^r$. Each configuration $p$ defines the
$n \times n$ matrix  $D_p=(d_{ij})$ = $(||p^i-p^j||^2)$, where
$||.||$ denotes the Euclidean norm. $D_p$ is called
the {\em Euclidean distance matrix (EDM)}
generated by configuration $p$. Obviously,
$D_p$ is a real symmetric matrix whose diagonal entries are all zeros.
A fundamental problem in distance geometry is to find out whether or not,
a given proper subset of the entries of $D_p$, the EDM generated by
configuration $p$,
suffices to uniquely determine the entire matrix $D_p$; i.e., to uniquely
recover $p$, up to a rigid motion. This problem
is known as the universal rigidity problem of bar frameworks.

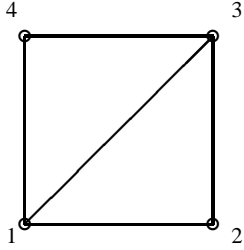
\begin{figure}
\thicklines
\setlength{\unitlength}{0.5mm}
\begin{picture}(100,70)(-100,-65)
\put(-30,-50){\circle{3}}
\put(20,-50){\circle{3}}
\put(-30,0){\circle{3}}
\put(20,0){\circle{3}}
\put(-30,-50){\line(1,0){50}}
\put(-30,-50){\line(0,1){50}}
\put(20,-50){\line(0,1){50}}
\put(-30,0){\line(1,0){50}}
\put(-30,-50){\line(1,1){50}}
\put(-35,-55){$1$}
\put(25,-55){$2$}
\put(25,5){$3$}
\put(-35,5){$4$}
\end{picture}
\caption{A bar framework $G(p)$ on 4 vertices in $\Rs^2$, where
         $V(G)=\{ 1,2,3,4\}$,
         $E(G)=\{ (1,2),(2,3),(3,4),(4,1),(1,3)\}$ and
        $p^1, p^2, p^3 , p^4$ are the vertices of the unit square.}
\label{e1}
\end{figure}

A {\em bar framework}, or framework for short, denoted by $G(p)$, in $\Rs^r$ is a
configuration $p$ in $\Rs^r$ together with a simple graph $G$
on the vertices $1,2, \ldots,n$. To avoid trivialities,
we assume throughout this chapter that graph $G$ is connected and not complete.
It is useful to think of each node $i$
of $G$ in a framework $G(p)$ as a universal joint located at $p^i$,
and of each edge $(i,j)$ of $G$ as a stiff bar of length $||p^i - p^j||$.
Hence, a bar framework is often defined as a collection of stiff bars
joined at their ends by universal joints.
Figure \ref{e1} depicts a framework $G(p)$ on 4 vertices in $\Rs^2$, where
$G$ is the complete graph $K_4$ minus an edge, and the points $p^1, \ldots,p^4$
are the vertices of the unit square.

We say that two frameworks $G(p)$ and $G(q)$
in $\Rs^r$ are {\em congruent} if $D_p=D_q$. Furthermore,
let $H$ denote the adjacency matrix of graph $G$, then
two frameworks $G(p)$ in $\Rs^r$ and $G(q)$ in $\Rs^s$
are said to be {\em equivalent} if $H \circ D_p= H \circ D_q$, where
$\circ$ denotes the {\em Hadamard product}, i.e., the element-wise product.
We say that framework $G(q)$ in $\Rs^r$ is {\em affinely-equivalent }
to framework $G(p)$ in $\Rs^r$ if $G(q)$ is equivalent to $G(p)$ and
configuration $q$ is obtained from configuration $p$ by an affine motion;
i.e., $q^i = A p^i + b$, for all $i=1,\ldots,n$,
for some $r \times r$ matrix $A$ and an $r$-vector $b$.

A framework $G(p)$ in $\Rs^r$ is said to be {\em universally rigid}
if every framework $G(q)$ in any dimension that
is equivalent to $G(p)$, is in fact congruent to $G(p)$; i.e., if
for every framework $G(q)$ in any dimension such that
$H \circ D_q= H \circ D_p$, it follows that  $D_q =  D_p$.

Thus, given $D_p =(d_{ij})$, the EDM generated by configuration $p$,
let $K \subset \{ (i,j) : i < j; \mbox{ for } i,j= 1,2,\ldots,n\}$.
Then the proper subset of entries of $D_p$ given by $\{d_{ij} : (i,j) \in K \}$
suffices to uniquely determine the entire matrix $D_p$ if and only if
framework $G(p)$ is universally rigid, where
$G=(V,E)$ is the graph with vertex set $V =\{1,2,\ldots,n\}$ and edge set $E= K $.
For example, the framework given in
Figure \ref{e1} is not universally rigid; and the subset of entries
of $D_p$ given by $\{d_{ij}: (i,j) \in E(G)\}$ does not uniquely
determine the entire matrix $D_p$ since the entry $d_{24}$ can
assume any value between $0$ and $2$.

The notion of dimensional rigidity is closely related to that of
universal rigidity.
A framework $G(p)$ in $\Rs^r$ is said to be {\em dimensionally rigid}
if there does not exist a framework $G(q)$ that is equivalent to $G(p)$,
in any Euclidean space of dimension $\geq r+1$. For example, the framework
$G(p)$ given in Figure \ref{e1} is obviously not dimensionally rigid since there
is an infinite number of frameworks $G(q)$ in $\Rs^3$ that are
equivalent to $G(p)$.

In this chapter, we survey some recently obtained results concerning
framework universal as well as dimensional rigidity.
These results are given in Section \ref{re} and their proofs
are given in Section \ref{spr}. Section \ref{prem} is dedicated to
the mathematical preliminaries needed for our proofs.
Our EDM approach of universal rigidity of bar frameworks extends to
the closely related notion of ``local" rigidity.
However, due to space limitation, ``local" rigidity
\cite{alf08} will not be considered here. Also, we will not consider
the other closely related notion of global rigidity \cite{con05,ght07}.

\section{Main Results}
\label{re}

The following theorem characterizes universal rigidity in terms
of dimensional rigidity and affine-equivalence.

\begin{theorem} [Alfakih \cite{alf07a}] \label{dimrig}
Let $G(p)$ be a bar framework on $n$ vertices in $\Rs^r$, $r \leq n-2$. Then
$G(p)$ is universally rigid if and only if the following two
conditions hold:
\begin{enumerate}
\item $G(p)$ is dimensionally rigid.
\item There does not exist a bar framework $G(q)$ in $\Rs^r$ that is
affinely-equivalent, but not congruent, to $G(p)$.
\end{enumerate}
\end{theorem}

The proof of Theorem \ref{dimrig} is given in Section \ref{spr}.
The notion of a stress matrix $S$ of a framework $G(p)$ plays an
important role in the characterization of universal rigidity of $G(p)$.
Let $G(p)$ be a framework on $n$ vertices in $\Rs^r$, $r \leq n-2$.
An {\em equilibrium stress} of $G(p)$ is a real valued function
$\omega$ on $E(G)$, the set of edges of $G$, such that
\beq \label{seq}
\sum_{j:(i,j) \in E(G)} \omega_{ij} (p^i - p^j) = \bz \mbox{ for all } i=1,\ldots,n.
\eeq

Let $\omega$ be an equilibrium stress of $G(p)$. Then the
$n \times n$ symmetric matrix $S=(s_{ij})$ where
\beq \label{defS}
s_{ij} = \left\{ \begin{array}{ll} -\omega_{ij} & \mbox{if } (i,j) \in E(G), \\
                        0   & \mbox{if } i \neq j \mbox{ and } (i,j) \not \in E(G), \\
                   {\displaystyle \sum_{k:(i,k) \in E(G)} \omega_{ik}} & \mbox{if } i=j,
                   \end{array} \right.
\eeq
is called the {\em stress matrix} associated with $\omega$, or a stress matrix
of $G(p)$.

Given framework $G(p)$ on $n$ vertices in $\Rs^r$, we define the following
$n \times r$ matrix
\beq \label{defP}
P := \left[ \begin{array}{c} {p^1}^T \\
                             {p^2}^T  \\
                              \vdots \\
                             {p^n}^T
         \end{array}  \right].
\eeq
$P$ is called the {\em configuration matrix} of $G(p)$.
Note that $P$ has full column rank since $p^1,\ldots,p^n$ affinely span $\Rs^r$.
The following lemma provides an upper bound on the rank of a stress matrix $S$.
\begin{lemma} \label{nullS}
Let $G(p)$ be a bar framework on $n$ nodes in $\Rs^r$, $r \leq n-2$,
and let $S$ and $P$ be a stress matrix and the configuration matrix
of $G(p)$ respectively. Then $SP=\bz$ and $Se = \bz $, where $e$ is
the vector of all 1's. Consequently, rank $S \leq n-r-1$.
\end{lemma}
\bpr It follows from (\ref{seq}) and (\ref{defS})
 that the $i$th row of $SP$ is given by
\[ s_{ii} {(p^i)}^T + \sum_{k=1, k \neq i }^n s_{ik} {(p^k)}^T =
\sum_{k:(i,k) \in E(G)} \omega_{ik} ({p^i} - {p^k})^T = \bz.
\]
Also, $e$ is obviously in the null space of $S$.
Hence, the result follows.
\epr

\subsection{Dimensional and Universal Rigidity In Terms of Stress  Matrices}

The following theorem provides a sufficient condition for the dimensional rigidity
of frameworks.

\begin{theorem} [Alfakih \cite{alf07a}] \label{suffdimrig}
Let $G(p)$ be a bar framework on $n$ vertices in $\Rs^r$ for some
$r\leq n-2$. If $G(p)$ admits a positive semidefinite stress matrix $S$
of rank $n-r-1$. Then $G(p)$ is dimensionally rigid.
\end{theorem}

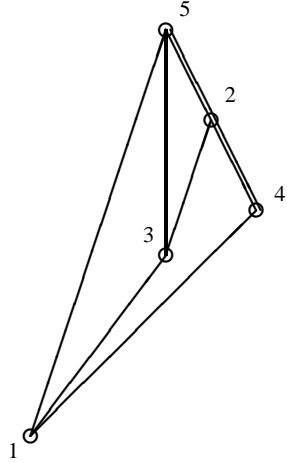
\begin{figure}
\thicklines
\setlength{\unitlength}{0.6mm}
\begin{picture}(140,110)(-100,-65)
\put(-30,-50){\circle{3}}
\put(10,20){\circle{3}}
\put(0,-10){\circle{3}}
\put(20,0){\circle{3}}
\put(0,40){\circle{3}}
\put(-30,-50){\line(1,3){30}}
\put(-30,-50){\line(3,4){30}}
\put(-30,-50){\line(1,1){50}}
\put(0,-10){\line(0,1){50}}
\put(0,-10){\line(1,3){10}}
\put(0,40){\line(1,-2){20}}
\put(1,40){\line(1,-2){20}}
\put(-35,-55){$1$}
\put(13,24){$2$}
\put(-5,-7){$3$}
\put(24,2){$4$}
\put(3,43){$5$}
\end{picture}
\caption{A dimensionally rigid framework $G(p)$ in $\Rs^2$
        (in fact $G(p)$ is also universally rigid)
        that does not admit a positive semidefinite stress matrix of rank 2.
        Note that the points $p^2,p^4$, and $p^5$ are collinear; i.e.,
        $G(p)$ is not in general position.}
\label{ce}
\end{figure}

The proof of Theorem \ref{suffdimrig} is given in Section \ref{spr}.
It is worth pointing out that the converse of Theorem \ref{suffdimrig} is not
true. Consider the following framework \cite{alf07a} $G(p)$
on 5 vertices in $\Rs^2$ (see Fig \ref{ce}),
where the configuration matrix $P$ is given by
\[ P = \left[ \begin{array}{rr} -3 & -5 \\ 1 & 2 \\ 0 & -1 \\
                      2 & 0 \\ 0 & 4 \end{array} \right], \]
and where the missing edges of $G$ are $(1,2)$ and $(3,4)$.
It is clear that $G(p)$ is dimensionally rigid (in fact $G(p)$ is also
universally rigid) while $G(p)$ has no positive
semidefinite stress matrix of rank 2.

The following result, which provides a sufficient condition for the
universal rigidity of a given framework, is a direct consequence of Theorems
\ref{dimrig} and \ref{suffdimrig}.

\begin{theorem}[Connelly \cite{con82,con99}, Alfakih \cite{alf07a}] \label{suff}
Let $G(p)$ be a bar framework on $n$ vertices in $\Rs^r$, for some $r \leq n-2$.
If the following two conditions hold:
\begin{enumerate}
\item $G(p)$ admits a positive semidefinite stress matrix $S$ of rank
$n-r-1$.
\item There does not exist a bar framework $G(q)$ in $\Rs^r$
that is affinely-equivalent, but not congruent, to $G(p)$.
\end{enumerate}
Then $G(p)$ is universally rigid.
\end{theorem}

A configuration $p$ (or a framework $G(p)$) is said to be {\em generic} if
all the coordinates of $p^1,\ldots,p^n$ are algebraically independent
over the integers. That is, if there does not exist
a non-zero polynomial $f$ with integer coefficients
such that $f(p^1,\ldots,p^n)=0$. Thus, for a generic framework, Theorem
\ref{suff} reduces to the following theorem.

\begin{theorem}[Connelly \cite{con99}, Alfakih \cite{alf10}] \label{suffg}
Let $G(p)$ be a generic bar framework on $n$ nodes in $\Rs^r$, for some $r \leq n-2$.
If $G(p)$ admits a positive semidefinite stress matrix $S$ of rank
$n-r-1$. Then $G(p)$ is universally rigid.
\end{theorem}

The proof of Theorem \ref{suffg} is given in Section \ref{spr}.
The converse of Theorem \ref{suffg} is also true.

\begin{theorem}[Gortler and Thurston \cite{gt09}] \label{neccg}
Let $G(p)$ be a generic bar framework on $n$ nodes in $\Rs^r$, for some $r \leq n-2$.
If $G(p)$ is universally rigid, then
there exists a positive semidefinite stress matrix $S$ of $G(p)$ of rank
$n-r-1$.
\end{theorem}

The proof of Theorem \ref{neccg} given in \cite{gt09} goes beyond the scope of
this chapter and will not be presented here.

At this point, one is tempted to ask whether a result similar to
Theorem \ref{suffg} holds if the genericity
assumption of $G(p)$ is replaced by the weaker assumption of general position.
A configuration $p$ (or a framework $G(p)$) in $\Rs^r$ is
said to be in {\em general position} if
no $r+1$ points in $p^1,\ldots,p^n$ are affinely dependent.
For example, a set of points in the plane are in general position
if no 3 of them are collinear. The following theorem answers this
question in the affirmative.

\begin{theorem}[Alfakih and Ye \cite{ay10}]  \label{suffgp}
Let $G(p)$ be a bar framework on $n$ nodes in general position in $\Rs^r$,
for some $r \leq n-2$.
If $G(p)$ admits a positive semidefinite stress matrix $S$ of rank
$n-r-1$. Then $G(p)$ is universally rigid.
\end{theorem}

The proof of Theorem \ref{suffgp} is given in Section \ref{spr}.
The following result shows that the converse of Theorem \ref{suffgp} holds for
frameworks $G(p)$ where graph $G$ is an $(r+1)$-lateration graph.
Such frameworks were shown to be universally rigid in \cite{zsy10}.
However, it is still an open question whether
the converse of Theorem \ref{suffgp} holds for frameworks
of general graphs.

A graph $G$ on $n$ vertices is
called an {\em $(r+1)$-lateration graph} if there is a permutation $\pi$ of the
vertices of $G$, $\pi(1), \pi(2), \ldots, \pi(n)$, such that
\begin{itemize}
   \item the first $(r+1)$ vertices, $\pi(1), \ldots, \pi(r+1)$, induce a clique
    in $G$, and
   \item each remaining vertex $\pi(j)$, for $j = (r+2),(r+3), \ldots, n$, is adjacent
                to $(r+1)$ vertices in the set $\{\pi(1), \pi(2), \ldots, \pi(j-1)\}$.
\end{itemize}

\begin{theorem}[Alfakih et al \cite{aty10c}] \label{lateral}
Let $G(p)$ be a bar framework on $n$ nodes in general position in $\Rs^r$,
for some $n \geq r+2$, where $G$ is an $(r+1)$-lateration graph. Then
there exists a positive semidefinite stress matrix $S$ of $G(p)$ of rank
$n-r-1$.
\end{theorem}

The proof of Theorem \ref{lateral} is given in Section \ref{spr}.
The preceding theorems have been stated in terms of stress matrices. The same
theorems can be equivalently stated in terms of Gale matrices,
as will be shown in the next subsection.

\subsection{Dimensional and Universal Rigidity in Terms of Gale  Matrices}

Let $G(p)$ be a framework on $n$ vertices in $\Rs^r$, $r \leq n-2$,
and let $P$ be the configuration matrix of $G(p)$.
Then the following $(r+1) \times n$ matrix
\beq \label{defPe}
{\cal P} := \left[ \begin{array}{c} P^T \\ e^T \end{array}  \right]
         = \left[ \begin{array}{ccc} p^1 &  \ldots & p^n \\
                           1   &  \ldots & 1
         \end{array}  \right]
\eeq
has full row rank since $p^1,\ldots,p^n$ affinely span $\Rs^r$.
Note that $r \leq n-1$. Let
\beq \label{defrr}
\rr = \mbox{the dimension of the null space
of } {\cal P}; \mbox{ i.e., } \rr=n-1-r.
\eeq
\begin{definition} \label{defLam}
Suppose that the null space of ${\cal P}$ is nontrivial, i.e., $\rr \geq 1$.
Any $n \times \rr$ matrix $Z$
whose columns form a basis of the null space of ${\cal P}$
is called a {\em Gale matrix} of configuration $p$ (or framework $G(p)$).
Furthermore, the $i$th row of $Z$, considered as
a vector in $\Rs^{\rr}$, is called a {\em Gale transform} of $p^i$ \cite{gal56}.
\end{definition}

Gale transform plays an important role in the theory of polytopes \cite{gru67}.
It follows from Lemma \ref{nullS} and (\ref{defS}) that
$S$ is a stress matrix of $G(p)$ if and only if
\beq \label{SZ1}
{\cal P} S = \bz, \mbox{ and } s_{ij}=0 \mbox{ for all }ij: i \neq j,(i,j) \not \in E(G).
\eeq
Equivalently,  $S$ is a stress matrix of $G(p)$ if and only if
there exists an $\rr \times \rr$ symmetric matrix $\Psi$ such that
\beq \label{SZ2}
S = Z \Psi Z^T , \mbox{ and }
s_{ij}= {(z^i)}^T \Psi z^j = 0 \mbox{ for all }ij: i \neq j,(i,j) \not \in E(G),
\eeq
where ${(z^i)}^T$ is the $i$th row of $Z$.
Therefore, the stress matrix
$S = Z \Psi Z^T$ attains its maximum rank of $\rr=n-1-r$ if and
only if $\Psi$ is nonsingular, i.e., rank $\Psi = \rr$, since
rank $S$ = rank $\Psi$.

Then Theorems \ref{suffdimrig}, \ref{suffg}, \ref{neccg} and \ref{suffgp}
can be stated in terms of Gale matrices as follows.

\begin{theorem} [Alfakih \cite{alf07a}]
Let $G(p)$ be a bar framework on $n$ vertices in $\Rs^r$ for some
$r\leq n-2$, and let $Z$ be a Gale matrix of $G(p)$. If there exists a positive definite
symmetric matrix $\Psi$ such that
\[ {(z^i)}^T \Psi z^j = 0 \mbox{ for each }ij: i \neq j, (i,j) \not \in E(G), \]
where ${(z^i)}^T$ is the $i$th row of $Z$.
Then $G(p)$ is dimensionally rigid.
\end{theorem}

\begin{theorem}
[Connelly \cite{con99}, Alfakih \cite{alf10}, Gortler and Thurston \cite{gt09}]
Let $G(p)$ be a generic bar framework on $n$ nodes in $\Rs^r$, for some $r \leq n-2$.
Let $Z$ be a Gale matrix of $G(p)$. Then
$G(p)$ is universally rigid if and only if there exists a positive definite
symmetric matrix $\Psi$ such that
\[ {(z^i)}^T \Psi z^j = 0 \mbox{ for each }ij: i \neq j, (i,j) \not \in E(G), \]
where ${(z^i)}^T$ is the $i$th row of $Z$.
\end{theorem}

\begin{theorem}[Alfakih and Ye \cite{ay10}]
Let $G(p)$ be a bar framework on $n$ nodes in general position in $\Rs^r$,
for some $r \leq n-2$. Let $Z$ be a Gale matrix of $G(p)$.
Then $G(p)$ is universally rigid if
there exists a positive definite symmetric matrix $\Psi$ such that
\[ {(z^i)}^T \Psi z^j = 0 \mbox{ for each }ij: i \neq j, (i,j) \not \in E(G), \]
where ${(z^i)}^T$ is the $i$th row of $Z$.
\end{theorem}

\section{Preliminaries}
\label{prem}

In this section we give the mathematical preliminaries needed for our proofs.
In particular, we review some basic terminology and results concerning
Euclidean distance matrices and affine motions. We begin with notation.

Throughout this chapter, $||.||$ denotes the Euclidean norm.
$|C|$ denotes the cardinality of a finite set $C$.
We denote the node set and the edge set of a simple graph $G$ by $V(G)$
and $E(G)$ respectively.
$\Sn$ denotes the space of $n \times n$ real symmetric matrices.
Positive semi-definiteness (positive definiteness) of a symmetric matrix
$A$ is denoted by $ A \succeq \bz$ ($A \succ \bz$).
For a matrix $A$ in $\Sn$, $\diag(A)$ denotes the $n$-vector formed
from the diagonal entries of $A$.
$e$ denotes the vector of all ones in $\Rs^n$.
$A \circ B$ denotes the Hadamard (element-wise) product of
matrices $A$ and $B$.
Finally, the $n \times n$ identity matrix is denoted by $I_n$;
and \bz denotes the zero matrix or the zero vector of the appropriate dimension.

\subsection{Euclidean Distance Matrices (EDMs)}

An $n \times n$ matrix $D=(d_{ij})$ is said to be a
{\em Euclidean distance matrix (EDM)} if and only if there exist points
$p^1, \ldots, p^n$ in some Euclidean space such that
$d_{ij}= \| p^i -p^j\|^2$ for all $i,j = 1, \ldots, n$.
The dimension of the affine subspace spanned by
$p^1, \ldots, p^n$ is called
the {\em embedding dimension} of $D$.

It is well known \cite{cri88,gow85,sch35,yh38} that a symmetric
$n \times n$ matrix $D$ whose diagonal entries are all zeros is EDM
if and only if $D$ is negative semidefinite on the subspace
\[ M := \{ x \in \Rs^{n} : e^T x = 0 \}, \]
where $e$ is the vector of all 1's.

Let $V$ be the $n \times (n-1)$ matrix whose columns form an orthonormal
basis of $M$; that is, $V$ satisfies:
\beq \label{defV} V^Te= \bz \;, \;\;\;\;  V^TV=I_{n-1} \;. \eeq
Then the orthogonal projection on $M$, denoted by $J$, is given by
$J:= VV^T = I_n - e e^T/n$.

Recall that $\snn$ denotes the subspace of symmetric matrices of order $n-1$ and
let $\Sh = \{ A \in \Sn : \diag(A) = \bz \}$. Consider the linear operator
$\TT_V: \Sh \rightarrow \snn$ such that

\beq \label{defTV}
\TT_V(D) := - \frac{1}{2} V^T D V,
\eeq
Then we have the following lemma.
\begin{lemma} [\cite{akw99}] \label{bb}
Let $D \in \Sh$. Then
$D$ is a Euclidean distance matrix of embedding dimension $r$
if and only if $\TT_V(D) \succeq \bz$ and rank $\TT_V(D) = r$.
\end{lemma}

Let $\KK_V: \snn \rightarrow \Sh$ defined by
\beq \label{defKV}
\KK_V(X):= \diag(VXV^T) \, e^T + e \, (\diag(VXV^T))^T - 2 \, VXV^T.
\eeq
Then it is not difficult to show that the operators $\TT_V$ and $\KK_V$
are mutually inverse \cite{akw99}. Thus, Lemma \ref{bb} implies that
$D$ in $\Sh$ is an EDM of embedding dimension $r$ if and only if
$D = \KK_V(X)$ for some positive semidefinite matrix $X$ of rank $r$.

Lemma \ref{bb} is used in the following subsection to characterize
the set of equivalent frameworks.

\subsection{Characterizing Equivalent Bar Frameworks}

Since all congruent frameworks have the same EDM, (or equivalently, the same
projected Gram matrix), in the rest of this chapter
we will identify congruent frameworks. Accordingly, for a given framework
$G(p)$ we assume without loss of generality that the centroid of the points
$p^1,\ldots,p^n$ coincides with the origin; i.e., $P^T e = \bz $, where
$P$ is the configuration matrix of $G(p)$.

Let $D=(d_{ij})$ be the EDM generated by
framework $G(p)$ in $\Rs^r$ and let $P$ be the configuration matrix
of $G(p)$ defined in (\ref{defP}).
Let $X = \TT_V(D)$, or equivalently, $D= \KK_V(X)$; and let $B=PP^T$
be the Gram matrix generated by the points $p^1,\ldots,p^n$.
Clearly, $B$ is positive semidefinite of rank $r$.
Observe that
\begin{eqnarray*}
d_{ij}= && || p^i -p^j||^2 , \\
      = && {(p^i)}^Tp^i+{(p^j)}^T p^j - 2 \; {(p^i)}^T p^j , \\
     = && (PP^T)_{ii} + (PP^T)_{jj} - 2  \; (PP^T)_{ij} .\\
\end{eqnarray*}
Therefore,
\[ D = \diag(B) e^T + e (\diag(B))^T - 2 B = \KK_V(X). \]
Hence,
\beq \label{defX}
B= V X V^T, \mbox{ and } X = V^T B V = V^T P P^T V .
\eeq
Furthermore, matrix $X$ is $(n-1) \times (n-1)$ positive semidefinite of rank $r$.
Accordingly, $X$ is called the {\em projected Gram matrix} of $G(p)$.

Now let $G(q)$ in $\Rs^s$ be a framework equivalent to $G(p)$. Let
$D_q$ and $D_p$ be the EDMs generated by $G(q)$ and $G(p)$ respectively.
Then $H \circ D_q = H \circ D_p$ where $H$ is the adjacency matrix of
graph $G$. Thus,
\beq \label{nullHKV}
H\circ (D_q-D_p)= H \circ \KK_V(X_q - X_p)= \bz,
\eeq
where $X_q$ and $X_p$ are the projected Gram matrices of $G(q)$ and $G(p)$
respectively.

Let $E^{ij}$ be the $n \times n$ symmetric matrix with 1's in the
$ij$th and $ji$th entries and zeros elsewhere. Further,
let
\beq \label{defM}
M^{ij} := \TT_V(E^{ij}) = - \frac{1}{2} V^T E^{ij} V .
\eeq
Then one can easily show that the set $\{M^{ij} : i \neq j, \; (i,j) \not \in E(G)\}$
forms a basis for the null space of $H \circ \KK_V$.
Hence, it follows from (\ref{nullHKV}) that
\beq \label{null2}
X_q - X_p = \sum_{ij: i \neq j, (i,j) \not \in E(G)} y_{ij} M^{ij},
\eeq
for some scalars $y_{ij}$. Therefore, given
a framework $G(p)$ in $\Rs^r$,
the set of projected Gram matrices of all frameworks $G(q)$ that
are equivalent to $G(p)$ is given by
\beq \label{defO}
\{ X :  X = X_p + \sum_{ij: i \neq j, (i,j) \not \in E(G)} y_{ij} M^{ij} \succeq \bz \}.
\eeq

The following lemma establishes the connection between Gale matrices and
projected Gram matrices.

\begin{lemma}[Alfakih \cite{alf01}] \label{vu}
Let $G(p)$ be a bar framework in $\Rs^r$ and let $P$ and $X$ be
the configuration matrix and the projected Gram matrix of $G(p)$ respectively.
Further, let $U$ and $W$ be the matrices whose columns form orthonormal bases
for the null space and the column space of $X$. Then
\begin{enumerate}
\item $VU$ is a Gale matrix of $G(p)$,
\item $VW = P Q$ for some $r \times r$ non-singular matrix $Q$.
\end{enumerate}
\end{lemma}
\bpr
It follows from (\ref{defX}) that $XU= V^T P P^T V U = \bz$. Thus $P^T VU =\bz$.
Hence, $VU$ is a Gale matrix of $G(p)$ since obviously $e^T VU =\bz$.

Now, $(VW)^T VU = \bz$. Thus $VW = PQ$ for some
matrix $Q$ since $P^T e= \bz$. Moreover, $Q$ is nonsingular since rank $PQ$ = $r$.
\epr

\subsection{Affine Motions}

Affine motions play an important role in the problem of universal rigidity
of bar frameworks. An affine motion in $\Rs^r$ is a map
$f: \Rs^r \rightarrow \Rs^r$ of the form
\[
f(p^i)= A p^i + b,
\]
for all $p^i$ in $\Rs^r$, where $A$ is
an $r \times r$ matrix and $b$ is an $r$-vector.
A rigid motion is an affine motion where matrix $A$ is orthogonal.

Vectors $v^1,\ldots,v^m$ in $\Rs^r$ are said to lie on a {\em quadratic
at infinity} if there exists a non-zero symmetric $r \times r$
matrix $\Phi$ such that
\beq
({v^i})^T \Phi v^i = 0 , \mbox{ for all } i=1,\ldots,m.
\eeq

The following lemma establishes the connection between the notion
of quadratic at infinity and affine motions.

\begin{lemma}(Connelly \cite{con05}) \label{conicp}
Let $G(p)$ be a bar framework on $n$ vertices in $\Rs^r$. Then the following
two conditions are equivalent:
\begin{enumerate}
\item There exists a bar framework $G(q)$ in $\Rs^r$ that is affinely-equivalent,
but not congruent, to $G(p)$,
\item The vectors $p^i-p^j$ for all $(i,j) \in E(G)$ lie on a quadratic at infinity.
\end{enumerate}
\end{lemma}

\bpr
Suppose that there exists a framework $G(q)$ in $\Rs^r$ that is affinely-equivalent,
but not congruent, to $G(p)$; and let
$q^i= A p^i + b$ for all $i=1,\ldots,n$.
Then $(q^i-q^j)^T(q^i-q^j)$ =  $(p^i-p^j)^TA^T A (p^i-p^j)$ =
$(p^i-p^j)^T(p^i-p^j)$ for all $(i,j) \in E(G)$.
Note that matrix $A$ is not orthogonal since $G(q)$ and $G(p)$ are not congruent.
Therefore, $(p^i-p^j)^T \Phi (p^i-p^j)$ = 0
for all $(i,j) \in E(G)$, where $\Phi= I_r- A^TA$ is a nonzero symmetric
matrix.

On the other hand, suppose that there exists a non-zero symmetric matrix $\Phi$
such that $({p^i-p^j})^T \Phi (p^i-p^j) = 0$, for all $(i,j) \in E(G)$.
Then $I_r- \delta \Phi \succ \bz $ for sufficiently small $\delta$.
Hence, there exists a matrix $A$ such that $I_r- \delta \Phi = A^T A$.
Note that matrix $A$ is not orthogonal since $\Phi$ is nonzero.
Thus, $({p^i-p^j})^T (I_r-A^T A) (p^i-p^j) = 0$ for all $(i,j) \in E(G)$.
Therefore, there exists a framework $G(q)$ in $\Rs^r$ that is equivalent
to $G(q)$, where $q^i= A p^i$ for all $i=1,\ldots,n$.
Furthermore, $G(q)$ is not congruent to $G(p)$ since $A$ is not orthogonal.
\epr

Note that Condition 2 in Lemma \ref{conicp} is expressed in terms of the
edges of $G$. An equivalent condition in terms of the missing edges of
$G$ can also be obtained using Gale matrices.
To this end,
let $\mm$ be the number of missing edges of graph $G$ and let $y=(y_{ij})$
be a vector in $\Rs^{\mm}$.
Let $\EE(y)$ be the $n \times n$ symmetric matrix whose $ij$th
entry is given by
\beq \label{defEE}
 \EE(y)_{ij} = \left\{ \begin{array}{cl}
                    y_{ij} & \mbox{ if } i \neq j \mbox{ and } (i,j) \not \in E(G),\\
                                   0 & \mbox{ Otherwise }.
    \end{array} \right.
\eeq
Then we have the following result.

\begin{lemma}(Alfakih \cite{alf10}) \label{lll}
Let $G(p)$ be a bar framework on $n$ vertices in $\Rs^r$ and let $Z$ be any
Gale matrix of $G(p)$. Then
the following two conditions are equivalent:
\begin{enumerate}
\item The vectors $p^i-p^j$ for all $(i,j) \in E(G)$ lie on a quadratic at infinity.
\item There exists a non-zero $y=(y_{ij}) \in \Rs^{\mm}$ such that:
\beq \label{conicZ}
V^T \EE(y) Z = \bz ,
\eeq
where $V$ is defined in (\ref{defV}).
\end{enumerate}
\end{lemma}

\bpr Let $P$ be the configuration matrix of $G(p)$, and let $U$ and $W$ be
the matrices whose columns form orthonormal bases for the null space and the
column space of $X$, the projected Gram matrix of $G(p)$. Then by
Lemma \ref{vu} we have
\begin{eqnarray*}
({p^i-p^j})^T \Phi (p^i-p^j) &&
= {(p^i)}^T \Phi p^i+{(p^j)}^T \Phi p^j-2 {(p^i)}^T \Phi p^j \\
 && = (P \Phi P^T)_{ii} +(P \Phi P^T)_{jj} -2 (P \Phi P^T)_{ij} \\
 && =  (VW \Phi' W^T V^T)_{ii} +(VW \Phi' W^T V^T)_{jj} -2 (VW \Phi' W^T V^T)_{ij} \\
 && = \KK_V(W \Phi' W^T )_{ij},
\end{eqnarray*}
where $\Phi'=Q \Phi Q^T$ for some nonsingular matrix $Q$, and where
$\KK_V$ is defined in (\ref{defKV}).

Therefore, $p^i-p^j$ for all $(i,j) \in E(G)$ lie on a quadratic at infinity
if and only if there exists a nonzero matrix $\Phi'$ such that
$H \circ \KK_V (W \Phi' W^T) = \bz$. But since the set
$\{M^{ij}: i \neq j,(i,j) \not \in E(G)\}$ forms a basis for the
null space of $H\circ \KK_V$, it follows that
vectors $p^i-p^j$ for all $(i,j) \in E(G)$ lie on a quadratic at infinity
if and only if there exists a nonzero $r \times r$ matrix $\Phi'$ and
a nonzero $y=(y_{ij})$ in $\Rs^{\mm}$ such that
\beq \label{nkv}
W \Phi' W^T = \sum_{ij:i\neq j,(i,j) \not \in E(G)} y_{ij} M^{ij}
    = - \frac{1}{2} \sum_{ij:i\neq j, (i,j) \not \in E(G)} y_{ij} V^T E^{ij} V
    = - \frac{1}{2} V^T \EE(y) V .
\eeq
Next we show that (\ref{nkv}) is equivalent to (\ref{conicZ}).
Suppose there exists a nonzero $y$ that satisfies (\ref{nkv}). Then by multiplying
(\ref{nkv}) from the right by $U$ we have that $y$ also satisfies (\ref{conicZ}).
Now suppose that there exists a nonzero $y$ that satisfies (\ref{conicZ}).
Then
\begin{eqnarray*}
V^T \EE(y) V & = [W \, U] \left[ \begin{array}{c} W^T \\ U^T \end{array} \right]
V^T \EE(y) V \; [W \, U] \left[ \begin{array}{c} W^T \\ U^T \end{array} \right], \\
 & = [W \, U] \left[ \begin{array}{cc} -2 \Phi' & 0 \\ 0 & 0  \end{array} \right]
\left[ \begin{array}{c} W^T \\ U^T \end{array} \right], \\
 & = -2 W \Phi' W^T.
\end{eqnarray*}
Thus $y$ also satisfies (\ref{nkv}) and the result follows.
\epr

\subsection{Miscellaneous Lemmas}

We conclude this section with the following lemmas that will be needed in our
proofs. We begin with the following well-known Farkas Lemma on the cone of positive
semidefinite matrices.

\begin{lemma} \label{farkash}
Let $A^1,\ldots,A^k$ be given $n \times n$ symmetric matrices. Then
exactly one of the following two statements hold:
\begin{enumerate}
\item there exists $Y \succ \bz$ such that trace $(A^i Y)=0$ for all $i=1,\ldots,k$.
\item there exists $x=(x_i) \in \Rs^k$ such that $x_1 A^1 + \cdots + x_k A^k \succeq \bz,
 \neq \bz$.
\end{enumerate}
\end{lemma}

\bpr Assume that statement 1 does not hold, and let
\[
 \LL = \{ Y \in \Sn : \trace (A^i Y) = 0, \mbox{ for all }i=1,\ldots,k \} .
\]
Then the subspace $\LL$ is disjoint from the interior of the cone of
$n \times n$ positive semidefinite matrices.
By the separation theorem \cite[page 96]{roc70}, there exists a nonzero
symmetric matrix $\Theta$ such that $\trace(\Theta Y) = 0$ for all $Y \in \LL$ and
$\trace(\Theta C) \geq 0$ for all $C \succ \bz$. Therefore, $\Theta \succeq \bz$ and
$\Theta = \sum_{i=1}^k \, x_i A^i $ for some nonzero $x =(x_i) \in \Rs^k$.
Hence, statement 2 holds.

Now assume that statements 1 and 2 hold and let
$\Theta= x_1 A^1 +\cdots+x_k A^k$. Then  on one hand, trace $(\Theta Y) > 0$;
and on the other hand trace ($\Theta Y)$ =
$ \sum_{i=1}^k x_i \trace (A^i Y)$ = 0, a contradiction.
Hence, the result follows.
\epr

The following lemma shows that Gale matrices have a useful property
under the general position assumption.

\begin{lemma} \label{Zli}
Let $G(p)$ be a bar framework on $n$ nodes in general position in $\Rs^r$ and
let $Z$ be any Gale matrix of $G(p)$.
Then any $\rr \times \rr$ sub-matrix of $Z$ is nonsingular.
\end{lemma}

\bpr
Assume $\rr \leq r$. The proof of the case where $\rr \geq r+1$
is similar. Let $Z'$ be any $\rr \times \rr$ sub-matrix of $Z$,
and without loss of generality, assume that it is the sub-matrix defined
by the rows $\rr+1,\rr+2,\ldots,2 \rr$. Then, $Z'$ is singular
if and only if there exists a nonzero $\xi \in \Re^{\rr}$ such
that $Z' \xi = \bz$. Clearly, $Z \xi $ is in the null space
of ${\cal P}$.
Furthermore, $Z' \xi = \bz$ if and only if the components
$(Z \xi)_{\rr+1}$ =  $(Z \xi)_{\rr+2}$ =  $\ldots $  =
$(Z \xi)_{2 \rr} =  \bz$. Now since $Z \xi \neq \bz$,
this last statement holds if and only if the following $r+1$ points
$p^1,p^2,\ldots,p^{\rr},p^{2 \rr+1},\ldots,p^{n}$ are affinely
dependent; i.e., $G(p)$ is not in general position.
\epr

\section{Proofs} \label{spr}
In this section we present the proofs of the theorems stated in Section \ref{re}.

\subsection{Proof of Theorem \ref{dimrig} }

Let $G(p)$ be a given framework on $n$ vertices in $\Rs^r$
for some $r\leq n-2$.
Clearly, if $G(p)$ is universally rigid then $G(p)$ is dimensionally rigid
and there does not exist a framework $G(p)$ in $\Rs^r$ that is affinely-equivalent,
but not congruent, to $G(p)$.

To prove the other direction, let $X_p$ be the projected Gram matrix of
$G(p)$. Let
$Q=[W \; U]$ be the orthogonal matrix whose columns are the eigenvectors
of $X_p$, where the columns of $U$ form an orthonormal basis for the null
space of $X_p$.

Now suppose that $G(p)$ is not universally rigid. Then there exists a framework
$G(q)$ in $\Rs^s$, that is equivalent, but not congruent, to $G(p)$,
for some $s$: $1 \leq s \leq n-1$.
Therefore, there exists a nonzero $\yy$ in $\Rs^{\mm}$ such that
$X(\yy)= X_p + \M(\yy) \succeq \bz$ where
$\M(\yy) = \sum_{(i,j) \not \in E} \yy_{ij} M^{ij}$.
Now for a sufficiently small positive scalar $\delta$ we have
\footnote{the rank function is lower semi-continuous on the set of
matrices of order $n-1$}.
\beq
X(t \yy) = X_p +  \M(t \yy)  \succeq \bz, \mbox{ and }
\rank (X(t \yy)) = \rank (X_p +  \M(t \yy)) \geq r,
\eeq
for all $t : 0 \leq t \leq \delta$. But,

\[
Q^T( X_p + \M(t \yy))Q = \left[ \begin{array}{rcc}
     \Lambda + t W^T  \M( \yy)  W  & & t W^T  \M(\yy) U \\
       t U^T  \M(\yy)  W & & t U^T  \M(\yy)  U
                                     \end{array} \right] \succeq \bz,
\]
where $\Lambda$ is the $r \times r$ diagonal matrix consisting
of the positive eigenvalues of $X_p$.
Thus $U^T \M(\yy) \, U \succeq \bz$ and
the null space of $U^T \M(\yy)\, U \subseteq$
the null space of $W^T \M(\yy) \, U$.

Therefore, if rank ($X(t_0 \yy)) \geq r+1$ for some $0 < t_0 \leq \delta$
we have a contradiction since $G(p)$ is dimensionally rigid. Hence,
rank ($X(t \yy)) = r$ for all  $t : 0 \leq t \leq \delta$. Thus,
both  matrices $U^T \M(\yy) U$ and  $W^T \M(\yy) U $ must be
zero. This implies that $\M(\yy)U = \bz$ i.e., $V^T \EE(\yy) Z= \bz$
which is also a contradiction by Lemma \ref{lll}.
Therefore, $G(p)$ is universally rigid.
\epr

\subsection{Proof of Theorem  \ref{suffdimrig} }

Let $G(p)$ be a given framework on $n$ vertices in $\Rs^r$
for some $r\leq n-2$ and let $Z$ be a Gale matrix of $G(p)$.
Let $X_p$ be the projected Gram matrix of $G(p)$, and let
$Q=[W \; U]$ be the orthogonal matrix whose columns are the eigenvectors
of $X_p$, where the columns of $U$ form an orthonormal basis for the null
space of $X_p$.

Assume that $G(p)$ admits a positive semidefinite stress matrix $S$
of rank $n-r-1$. Therefore, there exists a positive definite symmetric matrix
$\Psi$ such that ${(z^i)}^T \Psi z^j$ = 0  for all $ij: i \neq j, (i,j) \not \in E(G)$.
Hence, by lemma \ref{farkash}, there does not exist $y=(y_{ij}) \in \Rs^{\mm}$
such that
$\sum_{ij: i \neq j, (i,j) \not \in E(G)} y_{ij} (z^i{(z^j)}^T + z^j {(z^i)}^T)$
is a non zero positive semidefinite matrix. But
$z^i{(z^j)}^T + z^j {(z^i)}^T$ = $Z^T E^{ij} Z$. Thus,
there does not exist $y=(y_{ij}) \in \Rs^{\mm}$ such that $Z^T \EE(y)Z$ is a
nonzero positive semidefinite matrix. Hence,
there does not exist $y=(y_{ij}) \in \Rs^{\mm}$ such that $U^T \M(y) U $ is a
nonzero positive semidefinite matrix.

Now assume that $G(p)$ is not dimensionally rigid then there exists
a nonzero $y$ such that
$X = X_p +  \M(y)  \succeq \bz$ and rank $X \geq r+1$.
But
\[
Q^T( X_p + \M(y))Q = \left[ \begin{array}{rcc}
     \Lambda +  W^T \M( y) W  & &  W^T   \M(y) U \\
        U^T  \M(y)  W & &  U^T \M(y)  U
                                     \end{array} \right] \succeq \bz.
\]
Since $\Lambda + W^T \M(y) W$ is $r \times r$, it follows that
$U^T \M(y) U $ is a nonzero positive semidefinite, a contradiction.
\epr

\subsection{Proof of Theorem \ref{suffg} }

We begin with the following lemma.

\begin{lemma}(Connelly \cite{con05}) \label{conicpge}
Let $G(p)$ be a generic bar framework on $n$ vertices in $\Rs^r$. Assume that
each node of $G$ has degree at least $r$. Then
the vectors $p^i-p^j$ for all $(i,j) \in E(G)$ do not lie on a quadratic at infinity.
\end{lemma}

Now let $G(p)$ be a generic bar framework on $n$ vertices in $\Rs^r$. If
$G(p)$ admits a positive semidefinite stress matrix of rank $n-r-1$, then
each vertex of $G$ has degree at least $r+1$ \cite[Theorem 3.2]{alf07a}.
Thus Theorem \ref{suffg} follows from Lemmas \ref{conicp}
and \ref{conicpge} and Theorem \ref{suff}.

\subsection{Proof of Theorem \ref{suffgp} }

The main idea of the proof is to show that Condition 2 of Lemma \ref{lll}
does not hold under the assumptions of the theorem.
The choice of the particular Gale matrix to be used in equation (\ref{conicZ})
is critical in this regard.  The proof presented here is that given in \cite{ay10}.

Let $\bar{N}(i)$ denote the set of nodes of graph $G$ that are non-adjacent
to node $i$; i.e.,
\beq
\label{defNb}  \bar{N}(i) = \{ j \in V(G): j\neq i \mbox{ and } (i,j) \not \in E(G) \},
\eeq

\begin{lemma} \label{lemZh}
Let $G(p)$ be a bar framework on $n$ nodes in general position in $\Rs^r$,
$r \leq n-2$. Assume that $G(p)$ has a stress matrix $S$ of rank
$n-1-r $. Then there exists a Gale matrix $\hat{Z}$ of $G(p)$ such that
$\hat{z}_{ij}=0$ for all $j=1,\ldots,\rr$ and $i \in \bar{N}(j+r+1)$.
\end{lemma}

\bpr
Let $G(p)$ be in general position in $\Rs^r$ and assume that it has
a stress matrix $S$ of rank $\rr= (n-1-r )$.
Let $Z$ be any Gale matrix of $G(p)$, then it follows from
(\ref{SZ2}) that  $S = Z \Psi Z^T$ for some non-singular symmetric
$\rr \times \rr$ matrix $\Psi$.
Let us write $Z$ as:
\beq
Z = \left[ \begin{array}{c} Z_1 \\ Z_2 \end{array} \right],
\eeq
where $Z_2$ is $\rr \times \rr$. Then it follows from Lemma  \ref{Zli} that
$Z_2$ is non-singular. Now let
\beq \label{defZh}
\hat{Z}=(\hat{z}_{ij})= Z \Psi {Z_2}^T.
\eeq
Then $\hat{Z}$ is a Gale matrix of $G(p)$ since both $\Psi$ and $Z_2$
are non-singular. Furthermore,

\[
S = Z \Psi Z^T = Z \Psi \, [ Z_1^T \;\; Z_2^T] = [Z \Psi Z_1^T \;\; \hat{Z}].
\]
In other words, $\hat{Z}$ consists of the last $\rr$ columns of $S$.
Thus $\hat{z}_{ij} = s_{i,j+r+1}$. It follows by the definition of $S$ that
$s_{i,j+r+1}=0$ for all $i,j$ such that $i \neq (j+r+1)$ and
$(i,j+r+1) \not \in E(G)$. Therefore,
$\hat{z}_{ij} = 0$ for all $j=1,\ldots,\rr$ and $i \in \bar{N}(j+r+1)$.
\epr

\begin{lemma} \label{lemsys2}
Let the Gale matrix in (\ref{conicZ}) be $\hat{Z}$ as defined in (\ref{defZh}).
Then the system of equations (\ref{conicZ}) is equivalent to  the system
of equations
\beq \label{conicZ2}
\EE(y) \hat{Z} = \bz.
\eeq
\end{lemma}

\bpr
System of equations (\ref{conicZ}) is equivalent to the following system of
equations in the unknowns, $y_{ij}$ ($i \neq j$ and $(i,j) \not \in E(G)$)
and $\xi =(\xi_j) \in \Rs^{\rr}$:
\beq \label{eq1}
 \EE(y) \hat{Z} = e \, \xi^T.
\eeq
Now for $j=1,\ldots,\rr$, we have that the $(j+r+1,j)$th entry of
$\EE(y) \hat{Z}$ is equal to $\xi_j$. But using (\ref{defEE}) and Lemma \ref{lemZh}
we have
\[
(\EE(y) \hat{Z})_{j+r+1 ,j}  =  \sum_{i=1}^n \EE(y)_{j+r+1, i} \; \hat{z}_{ij}
                   = \sum_{i: i \in \bar{N}(j+r+1)} y_{j+r+1, i} \; \hat{z}_{ij}  =0.
\]
Thus, $\xi = 0$ and the result follows.
\epr

\begin{lemma} \label{Zaffine}
Let $G(p)$ be a bar framework on $n$ nodes in general position in $\Rs^r$,
$r \leq n-2$.
Assume that $G(p)$ has a positive semidefinite stress matrix $S$ of rank
$\rr= n-1-r $. Then there does not exist a framework $G(q)$ in $\Rs^r$ that
is affinely-equivalent, but not congruent, to $G(p)$.
\end{lemma}

\bpr Under the assumption of the lemma, we have that
deg$(i) \geq r+1$ for all $i \in V(G)$,
i.e., every node of $G$ is adjacent to at least $r+1$ nodes
(for a proof see \cite[Theorem 3.2]{alf07a}). Thus
\beq \label{deg}
|\bar{N}(i)| \leq n-r-2 = \rr - 1 \mbox{ for all } i \in V(G).
\eeq

Furthermore,
it follows from  Lemmas \ref{lemZh}, \ref{lemsys2}
and \ref{lll} that
the vectors $p^i-p^j$ for all $(i,j) \in E(G)$ lie on a quadratic at infinity
if and only if system of equations (\ref{conicZ2}) has a non-zero solution $y$.
But (\ref{conicZ2}) can be written as

\[
\sum_{j: \in \bar{N}(i)} y_{ij} \hat{z}^j = 0 , \mbox{ for } i=1,\ldots,n,
\]
where $(\hat{z}^i)^T$ is the $i$th row of $\hat{Z}$. Now  it follows from (\ref{deg})
that  $y_{ij}=0$ for all $(i,j) \not \in E(G)$ since by Lemma \ref{Zli}
any subset of $\{\hat{z}^1,\ldots,\hat{z}^n\}$ of cardinality $\leq \rr-1$ is linearly
independent.

Thus system (\ref{conicZ2}) does not have a nonzero solution $y$. Hence
the vectors $p^i-p^j$, for all $(i,j) \in E(G)$, do not lie on a quadratic at infinity.
Therefore, by Lemma \ref{conicp}, there does not exist a framework $G(q)$ in $\Rs^r$
that is affinely-equivalent, but not congruent, to $G(p)$.
\epr

Thus, Theorem \ref{suffgp} follows from
Lemma \ref{Zaffine} and Theorem \ref{suff}.

\subsection{Proof of Theorem \ref{lateral} }

The proof of Theorem \ref{lateral} is constructive, i.e.,
an algorithm is presented to construct the desired stress matrix.
The proof presented here is a slight modification
of that given in \cite{aty10c}.

Let $G(p)$ be a framework on $n$ vertices in general position
in $\Rs^r$, $n \geq r+2$, and let $Z$ be a Gale matrix of $G(p)$.
An $n \times n$ symmetric matrix $S$ that satisfies
\[
{\cal P} S = 0 , \mbox{ or equivalently } S=Z \Psi Z^T
   \mbox{ for some symmetric matrix } \Psi,
\]
is called a {\em pre-stress matrix},
where ${\cal P}$ is defined in (\ref{defPe}).
Thus, it follows from (\ref{SZ1}) and (\ref{SZ2}) that
$S$ is a stress matrix of $G(p)$ if and only if  $S$ is a pre-stress matrix
and $s_{ij}=0$ for all $ij: i \neq j$,  $(i,j) \not \in E(G)$.

Clearly, $S^n=ZZ^T$ is
a positive semidefinite pre-stress matrix of rank $\rr=n-r-1$.
If $S^n$ satisfies $s^n_{ij}=0$ for all $ij: i \neq j$, $(i,j) \not \in E(G)$,
then we are done since $S^n$ is the desired stress matrix.
Otherwise, if $S^n$ is not a stress matrix, we need to zero out the entries
which should be zero but are not, i.e., the entries $s^n_{ij} \neq 0$, $i \neq j $
and $(i,j)\not \in E(G)$.  We do this in reverse order by column (row); first, we
zero out the entries $s^n_{in} \neq 0$, for $i<n$ and $(i,n)\not \in E(G)$, and then
do the same for columns (rows) $(n-1), (n-2), \ldots, (r+3)$.  This ``purification''
process will keep the pre-stress matrix positive semidefinite
and maintain rank $n-r-1$.

Let $G$ be an $(r+1)$-lateration graph with lateration order $1,2,\ldots,n$;
i.e., the vertices, $1, 2, \ldots, r+2$, induce a clique in $G$, and
each remaining vertex $k$, for $k = r+3, \ldots, n$, is adjacent
to $(r+1)$ vertices in the set $\{1, 2, \ldots, k-1\}$.
Let
\beq \label{Np}
\bar{N}'(k)= \{ i \in V(G): i < k \mbox{ and } (i, k) \not \in E(G) \}.
\eeq
Then  for $k=r+3, \ldots, n$,
\beq
| \bar{N}'(k)| = k-r-2,
\eeq

We first show how to purify the last column (or row) of $S^n=ZZ^T$.
Let $Z^n$ denote the sub-matrix of $Z$ obtained by keeping only rows
with indices in $\bar{N}'(n) \cup \{n\}$. Then $Z^n$ is a square matrix of order
$\rr= n-r-1$. Furthermore, by Lemma \ref{Zli}, it follows that $Z^n$ is nonsingular.
Let $b^n$ denote the vector in $\Rs^{\rr}$ such that
\[  b^n_i = \left\{ \begin{array}{cl} -s^n_{in} & \mbox{if } i \in \bar{N}'(n),\\
                                      1    & \mbox{if } i=n.
              \end{array} \right.
\]
Now let $\xi_n \in \Rs^{\rr}$ be the unique solution of the system of equations
\[
Z^n \xi_n = b^n.
\]

\begin{lemma}
        \label{stepn}
Let $S^{n-1} = S^n + Z \, \xi_n {\xi_n}^T Z^T = Z (I + \xi_n {\xi_n}^T ) Z^T$.
Then
 \begin{enumerate}
  \item $S^{n-1}$ is a pre-stress matrix of $G(p)$, i.e., ${\cal P} S^{n-1}=0$.
  \item $S^{n-1}\succeq \bz$ and the rank of $S^{n-1}$ remains $n-r-1$.
  \item $s^{n-1}_{in}=0$ for all $i: i< n,\ (i,n)\not\in E(G)$.
\end{enumerate}
\end{lemma}
\bpr
The first statement is obvious.
The second statement follows since $I + \xi_n \xi_n^T \succ \bz$.
The third statement is also true by construction. For all
$i< n,\ (i,n)\not \in E(G)$, i.e., for all $i \in \bar{N}'(n)$, we have
$s^{n-1}_{in}= s^n_{in} + b^n_i b^n_n = s^n_{in} -  s^n_{in} = 0$.
\epr

We continue this purification process for columns $(n-1),\ldots,k, \ldots, (r+3)$.
Before the $k$th purification step, we have $S^{k}\succeq \bz$,
${\cal P} S^{k}=0$, rank $S^{k}  = n-r-1$, and
        \[s^{k}_{ij}=0, \mbox{ for all }ij: i \neq j , \; (i,j) \not \in E(G),
         \mbox{ and for all } j = k+1,\ldots,n . \]

Let $Z^k$ denote the sub-matrix of $Z$ obtained by keeping only rows
with indices in $\bar{N}'(k) \cup \{k, k+1, \ldots, n \}$.
Then $Z^k$ is a square matrix of order $\rr= n-r-1$.
Furthermore, by Lemma \ref{Zli}, it follows that $Z^k$ is nonsingular.
Let $b^k$ denote the vector in $\Rs^{\rr}$ such that
\[  b^k_i = \left\{ \begin{array}{cl} -s^{k}_{ik} & \mbox{if } i \in \bar{N}'(k),\\
                                      1    & \mbox{if } i=k, \\
                                      0    & \mbox{if } i=k+1,\ldots,n.
              \end{array} \right.
\]
Now let $\xi_k \in \Rs^{\rr}$ be the unique solution of the system of equations
\[
Z^k \xi_k = b^k.
\]

The following lemma shows results analogous to those in Lemma
\ref{stepn}, for the remaining columns.
\begin{lemma} \label{stepk}
        Let $S^{k-1}=S^k+ Z \, \xi_k \xi_k^T Z^T$. Then
 \begin{enumerate}
  \item $S^{k-1}$ is a pre-stress matrix of $G(p)$, i.e., ${\cal P} S^{k-1}=0$.
  \item $S^{k-1}\succeq \bz$ and the rank of $S^{k-1}$ remains $n-r-1$.
  \item $s^{k-1}_{ij}=0$ for all $i: i< j, \; (i,j) \not \in E(G)$ and for
                      all $j=k,\ldots,n$.
\end{enumerate}
\end{lemma}
\bpr
The proof of the first two statements is identical to that in Lemma \ref{stepn}.
The third statement is again true by construction.
For each $i< k,\ (i,k)\not \in E(G)$, i.e., for all $i \in \bar{N}'(k)$,
we have
\[
s^{k-1}_{ik}= s^k_{ik} + b^k_i b^k_k = s^k_{ik} -  s^k_{ik} = 0.
\]
Furthermore, for $j = k+1, \ldots, n$, the $j$th
column (or row) of $Z \, \xi_k \xi_k^T Z^T$ has all zero entries,
which means that the entries
in the $j$th column (or row) of $S^{k-1}$ remain unchanged from $S^k$.
\epr

\noindent {\bf Proof of Theorem \ref{lateral} }

The matrix
\[ S^{r+2} = S^{r+3}+Z \, \xi_{r+3} \xi_{r+3}^T Z^T =
  Z (I + \xi_n \xi_n^T + \cdots + \xi_{r+3} \xi_{r+3}^T) Z^T,
\]
obtained at the ``$(r+3)$th" step of the above process, is by  Lemmas
\ref{stepn} and \ref{stepk} a positive semidefinite pre-stress matrix
of rank $n-r-1$.
Furthermore, $s^{r+2}_{ij} = 0 $ for all $ij: i \neq j, \; (i,j) \not \in E(G)$ and
for all $j=r+3, r+4,\ldots,n$. But since the vertices $1,2,\ldots,r+2$ induce
a clique in $G$, it follows that
\[ s^{r+2}_{ij} = 0 \mbox{ for all } ij: i \neq j, \; (i,j) \not \in E(G).  \]
Hence, $S = S^{r+2}$ is a positive semidefinite stress matrix of $G(p)$ of
rank $n-r-1$; i.e., $S^{r+2}$ is the desired stress matrix.

\begin{acknowledgement}
Research supported by the Natural Sciences and Engineering
 Research Council of Canada.
\end{acknowledgement}

\end{document}